\documentclass[11pt]{amsart}
\usepackage{amssymb,latexsym,amsfonts,amsthm,amsmath}
\usepackage{enumerate}
\usepackage{epic, curves}

\numberwithin{equation}{section}

\theoremstyle{definition}
\newtheorem{example}[equation]{Example}
\newtheorem{definition}[equation]{Definition}

\theoremstyle{plain}
\newtheorem{corollary}[equation]{Corollary}
\newtheorem{lemma}[equation]{Lemma}

\newtheorem{theorem}[equation]{Theorem}

\theoremstyle{remark}
\newtheorem{remark}[equation]{Remark}

\begin{document}
\title{Aleksandrov Surfaces and Hyperbolicity}
\author{Byung-Geun Oh}
\thanks{Supported by NSF Grant DMS-0244421}
\address{Department of Mathematics, Purdue University, West Lafayette, IN 47907}
\email{boh@math.purdue.edu}
\date{February 18, 2004}

\begin{abstract}
Aleksandrov surfaces are a generalization of two dimensional Riemannian manifolds, and it
is known that every open simply connected Aleksandrov surface is conformally
equivalent either to the unit disc (hyperbolic case) or to the plane (parabolic
case). We prove a criterion for hyperbolicity of Aleksandrov surfaces which
have \emph{nice} tilings and where negative curvature dominates. We then apply
this to generalize a result of Nevanlinna and give a partial answer for his
conjecture about line complexes.
\end{abstract}
\maketitle

\newcommand{\D}{\Delta}
\newcommand{\wD}{\widetilde{\Delta}}
\newcommand{\T}{\mathcal{T}}
\newcommand{\PP}{\mathcal{P}}

\section{Introduction}\label{S:Intro}
An \emph{Aleksandrov surface} $S$ is an open simply-connected two-dimensional topological manifold
equipped with an intrinsic metric\footnote{A metric is called \emph{intrinsic} if the distance between two 
points in the space is the infimum of the lengths of curves connecting 
these points.} whose length element is locally defined by the form
\begin{equation}\label{eqre}
 e^{u(z)}|dz|, 
\end{equation}
where $z$ is a complex local coordinate and $u$ is a difference of two subharmonic
functions such that $\exp u$ is locally integrable on rectifiable curves in the $z$-plane.  
A local coordinate for which \eqref{eqre} holds is called an \emph{isothermal} coordinate,
and isothermal coordinates define a complex analytic structure on every Aleksandrov
surface. In fact by the Uniformization Theorem of Huber \cite{Hu}, there is an isometry 
\begin{equation}\label{hu}
h: S \to D(R) := \{z \in \mathbb{C} : |z| < R \},\; \;  R \in (0, \infty],
\end{equation}
where $D(R)$ is equipped with a length element of the form \eqref{eqre}. Therefore one
can study complex analysis on Aleksandrov surfaces, and in this paper we give a 
criterion for hyperbolicity (the case $R < \infty$ in \eqref{hu}).

To state our main theorem, we need to introduce several definitions. First,
the \emph{integral curvature} $\omega$ of $S$ is a signed Borel measure defined by the
negative of the generalized Laplacian of $u$ (i.e.,  $\omega = -\D u$), where $u$ is as
in \eqref{eqre}. The \emph{total angle} $T(v)$ at a point $v \in S$ is 
\begin{equation}\label{ta}
   T(v) = 2\pi - \omega(\{v\}),
\end{equation}
and we define the \emph{angular curvature} $K(\D)$ of a triangle
$\D$ with 
vertices $v_i$, $i=1,2,3$, as the quantity
\begin{equation}\label{ac}
K(\D) = 2\pi \sum_{i=1}^{3} \frac{\theta(v_i, \D)}{T(v_i)} - \pi,
\end{equation}
where $\theta(v_i, \D)$ is the angle of $\D$ at $v_i$.
A tiling $\T$ of $S$ is a set of triangles in $S$ which covers $S$ locally finitely and such that
for different triangles $\D$ and $\D'$ in $\T$, $\D \cap \D'$ 
is either empty, a set of common vertices, or a set of common sides.
Finally, we define a \emph{partition} $\PP$ of $\T$ as follows: each element
$C_{\alpha} \in \PP$, called a \emph{cluster},  is a finite union of 
triangles in $\T$ with connected 
interior such that $C_{\alpha} \cap C_{\beta}$ contains no triangle in $\T$
when $\alpha \ne \beta$ and $\bigcup_{C_{\alpha} \in \PP} C_{\alpha} 
= S$. Let $\# (C_\alpha)$ be the number of triangles in $C_\alpha$.

\begin{theorem}\label{T}
 Suppose that an Aleksandrov surface $S$ is of curvature at most 
 $k \in \mathbb{R}$.\footnote{This concept is defined in Section~\ref{S:Fact}, 
 and if $S$ is a smooth surface then it is equivalent to the boundedness
 of the Gaussian curvature from above by $k$.}
 If there are a tiling $\T$ and a partition $\PP$ on $S$ such that 
 for some constants $\epsilon > 0$ and $M > 0$,
 \begin{align*}
  &\mbox{\emph{(M1) }} \# (C_\alpha) \leq M \mbox{ for every }C_{\alpha} 
     \in \PP, \\
  &\mbox{\emph{(M2) }} \sum_{\D \subseteq C_{\alpha}} K(\D) 
     \leq -\epsilon \pi \mbox{ for every }C_{\alpha} \in \PP, \\
  &\mbox{\emph{(R1) }} \theta(v, \D) \geq \epsilon \mbox{ for every }
      \D \in \T \mbox{ and every vertex } v \mbox{ of } \D,   \\
  &\mbox{\emph{(R2)} if } k > 0, \mbox{ then \emph{(the perimeter of $\D$)}} 
      \leq (2\pi - \epsilon)k^{-1/2} \mbox{ for every }\D \in \T,
 \end{align*}
 then $S$ is hyperbolic; i.e., conformally equivalent to the unit disc.
\end{theorem}

The meaning of our main conditions (M1) and (M2) is that negative curvature dominates 
$S$ uniformly, and note that $S$ is allowed to  have positive integral curvature
in some subsets.  

One cannot drop (R1), one of the regularity conditions, as Example \ref{Ex} below shows.

\begin{example}\label{Ex}
 For every $\epsilon > 0$ there exists a parabolic Aleksandrov surface of curvature at most zero which 
 has a tiling $\T$ such that every triangle in $\T$ is isometric to a Euclidean triangle of area at most $\epsilon$ and 
 $T(v) = 4\pi$ at every vertex $v$ of triangles in $\T$.
\end{example}

With $M=1$, it is easy to see that 
the surface described in Example~\ref{Ex} 
satisfies all assumptions of Theorem~\ref{T} except (R1).
This example also shows that in the
Proposition~1.7 of \cite{BE2} the assumption about diameters of 
triangles cannot be replaced by an assumption about areas, while our 
Theorem~\ref{T} implies that it can be replaced by an assumption about angles.
Example~\ref{Ex} will be constructed in Section~\ref{S:Ex}.

When $k > 0$, we don't know whether we can drop the condition (R1) or (R2). However, an Aleksandrov surface of curvature
at most $k>0$ may be parabolic if it satisfies only the main conditions (M1) and (M2)
of Theorem~\ref{T}, as Example~1.2 of \cite{BE1} shows: in fact,
the surface constructed in \cite{BE1} is a surface over $\mathbb{C}$ (with
the induced Euclidean metric), but one can also realize it as a surface over
$\overline{\mathbb{C}}$ (with the induced spherical metric, which is not 
complete). Then this surface becomes an Aleksandrov surface of
curvature at most 1, and one can check that it satisfies the
conditions (M1) and (M2) only.

Let $X$ be a Riemann surface and $f: X \to \overline{\mathbb{C}}$ a non-constant
holomorphic function. For each $p \in X$, we define $d_f(p)$ as the 
radius of the largest open spherical disc centered at $f(p)$
where the inverse of $f$ exists.
Then the \emph{spherical Bloch constant} $\mathcal{B}_f$ of $f$ is defined by
the formula
$\mathcal{B}_f := \sup_{p \in X} d_f(p)$.  

\begin{theorem}[M.~Bonk and A.~Eremenko \cite{BE1}, Proposition 8.4]\label{T:BE}
 Let X be a Riemann surface and $f: X \to \overline{\mathbb{C}}$ a non-constant holomorphic
function without asymptotic values such that $\mathcal{B}_f \leq \pi/2 - \epsilon$
for some $\epsilon > 0$. Then there exists a tiling $\mathcal{T}$ of $X$ with respect to
the pull-back spherical metric such that the image $f(\D)$ of every triangle $\D \in \mathcal{T}$ is contained
in a closed spherical disc of radius $\mathcal{B}_f$ and the set of all the vertices of the triangles in $\mathcal{T}$
is equal to the set of critical points of $f$.
\end{theorem}

Theorem~\ref{T:BE} gives a sufficient condition for the existence of a tiling, which is \emph{apriori} assumed in 
our Theorem~\ref{T}. Also note that if $S$ is a Riemann surface spread over $\overline{\mathbb{C}}$ and 
ramified only over finitely many points, then a tiling always exists with respect to the pull-back spherical
metric (cf. Section~\ref{S:Appl}).

M. Bonk and A. Eremenko conjectured in \cite{BE2} (cf. \cite{Bo}) the following statement:

\medskip
\noindent \emph{Suppose that $S$ is an open simply-connected complete Aleksandrov surface with a tiling $\mathcal{T}$ such that
for some $q \in (1,3]$ and $\epsilon > 0$, the following two conditions are satisfied:\\
\noindent (1) every triangle $\D \in \mathcal{T}$ is isometric to a spherical triangle with circumscribed radius
at most 
\begin{equation}\label{strange}
  \mathcal{R}_{q,\epsilon} :=\arctan \sqrt{\frac{-\cos (\pi q / 2)}{\cos^3 (\pi q / 6)}} - \epsilon ;
\end{equation}
(2) the total angle at each vertex is at least $2\pi q$. Then $S$ is of hyperbolic type.}
\medskip
 
\noindent The expression $\mathcal{R}_{q,\epsilon}$ in 
\eqref{strange} with $\epsilon=0$ is the circumscribed radius of a spherical
equilateral triangle with angles $\pi q / 3$.

M. Bonk and A. Eremenko proved this conjecture in \cite{BE2} for the cases $q=2, 3$
and for the limiting case $q \to 1$ with an appropriate interpretation.  
This conjecture still remains open for general $q \in (1,3]$, but
as one of the applications of Theorem~\ref{T} we 
have the following corollary:

\begin{corollary}
Suppose $S$ is an Aleksandrov surface as described in the above conjecture with an additional assumption
\emph{(R1)} of Theorem~\ref{T}. Then $S$ is of hyperbolic type.
\end{corollary}

\begin{proof}
It suffices to check the four conditions in Theorem~\ref{T} because every triangle in $\T$ is spherical.
In fact, we can guarantee that $S$ is of curvature at most 1 only for $q \in [2,3]$, but Theorem~\ref{T} can be applied 
to all $q \in (1,3]$ (with $k=1$) because singular points\footnote{A point $p \in S$ is called \emph{singular} 
if $T(p) \neq 2 \pi$.} of $S$ are contained in the set of vertices of triangles in $\T$. 
See Remark~\ref{rem} for details.

First, we observe that the area of every spherical triangle with circumscribed radius at most
$\mathcal{R}_{q,\epsilon}$ is less than or equal to $\pi (q-1) - \eta$, where
$\eta > 0$ is a constant depending only on $\epsilon$. Therefore every triangle $\D \in \T$ has
angular curvature 
\[
K(\D) \leq 2 \pi \frac{ \pi q - \eta}{2 \pi q} - \pi = -\frac{\eta}{q},
\]
so the conditions (M1) and (M2) of Theorem~\ref{T} are satisfied with $M=1$. Note that we already assumed 
the condition (R1) in this corollary. The condition (R2) is trivially satisfied 
because the circumscribed radius of every triangle $\D \in \T$ is at most $\mathcal{R}_{q,\epsilon} \leq \pi/2 - \epsilon$.
\end{proof}

More applications will be given in Section~\ref{S:Appl}.

\section{Preliminaries}
\subsection{Aleksandrov surfaces of curvature at most $k$}\label{S:Fact}
Let $S$ be an Aleksandrov surface with
an intrinsic metric $\rho$. A curve joining points $x, y \in S$ is called a
\emph{shortest curve} if its length is equal to the distance $d_\rho (x,y)$
between $x$ and $y$, and is denoted by $[x,y]$. 
(This notation is ambiguous because we do not require the uniqueness of shortest curves.)

A \emph{geodesic triangle} $\D = \D(v_1 v_2 v_3)$ in $S$ is a closed set 
which is homeomorphic to the closed unit disc and whose boundary $\partial \D$ 
consists of three shortest curves $[v_1, v_2], [v_2, v_3],[v_3, v_1]$, called
\emph{sides} of $\D$. Each point $v_i$ is called a \emph{vertex} of $\D$,
and the \emph{perimeter} of $\D$ is the sum 
$d_\rho (v_1, v_2) + d_\rho (v_2, v_3) + d_\rho (v_3, v_1)$.
A geodesic triangle $\D$ is called a \emph{simple $($geodesic$)$ triangle}
if for any two points $x, y \in \D$ there exists a shortest curve
$[x,y] \subseteq \D$. 

Now suppose $L$ and $M$ are two curves in $S$ that have a common starting 
point $o$. On $L$ and $M$ we choose arbitrary points $x$ and $y$ respectively
and let 
\[
X = d_\rho(o,x), \quad Y=d_\rho(o,y), \quad Z=d_\rho(x,y).
\]
We then construct a triangle in the Euclidean plane with side-lengths
$X, Y$ and $Z$, and let $\gamma_{LM}(X,Y)$ be the angle opposite to the
side of length $Z$ in this Euclidean triangle. Then the \emph{upper angle}
between $L$ and $M$ is 
\begin{equation}\label{upperangle}
 \overline{\angle (L,M)} = \limsup_{X,Y \to 0} \gamma_{LM}(X,Y) \in [0, \pi].
\end{equation}
Using this definition of upper angles, one can define (upper) angles of
a geodesic triangle at each of its vertices. We make a remark that 
we will mostly use the terminology ``angles'' instead of ``upper angles''
when they are associated to ``small'' simple triangles, because for any two
``small'' shortest curves the angle between them exists; i.e., the limit
of $\gamma_{LM}(X,Y)$ in \eqref{upperangle} exists as $X$ and $Y$ go to zero
(\cite{BN}, pp. 176--177, 182, 187).

Now fix a number $k \in \mathbb{R}$. We denote by $S(k)$ the Euclidean plane
when $k=0$, the hyperbolic plane of curvature $k$ when $k<0$, and a two
dimensional open hemisphere of radius $k^{-1/2}$ when $k>0$. Then for a given
geodesic triangle $\D$ in $S$, we associate a new triangle $\D_k$ on $S(k)$
which has the same side-lengths as $\D$. (When $k \leq 0$ the 
triangle $\D_k$ exists on account of the triangle inequality, and when
$k>0$ it is necessary to require that the perimeter of $\D$ is less than
$2 \pi k^{-1/2}$. When $k>0$ it will be assumed throughout that
this requirement is satisfied.) If $\D$ has upper angles $\overline{\alpha},
\overline{\beta}, \overline{\gamma}$ and the corresponding angles of $\D_k$
are $\alpha_k, \beta_k, \gamma_k$, then the \emph{$k$-excess} of $\D$ is the 
quantity
\[
\delta_k(\D) = (\overline{\alpha} + \overline{\beta} + \overline{\gamma})
 - (\alpha_k + \beta_k + \gamma_k).
\]

\begin{definition}\label{domainrk}
A domain $R_k$ is an intrinsic metric space satisfying the following properties: 
\begin{enumerate}[(a)]
 \item any two points in $R_k$ can be joined by a shortest curve;
 \item each geodesic triangle in $R_k$ has non-positive $k$-excess;
 \item if $k>0$, the perimeter of each geodesic triangle in $R_k$ 
       is less than $2 \pi k^{-1/2}$.
\end{enumerate}
\end{definition}

We say that an Aleksandrov surface $S$ is of curvature at most $k$
if every point in $S$ has a neighborhood which is a domain $R_k$, 
and denote by $\mathfrak{S}(k)$ the collection 
of Aleksandrov surfaces of curvature at most $k$.
Note that this notion is equivalent to the boundedness of
the Gaussian curvature from above (by $k$) when $S$ is a Riemannian surface,
and also note that if an Aleksandrov surface has a point $v$ such that
$\omega (\{ p \}) >0$, or equivalently $T(p) < 2 \pi$,  then it
does not belong to $\mathfrak{S}(k)$ for any $k < \infty$.

\subsection{Generalized simple geodesic triangles}\label{S:gentri}
Let $\overline{S}$ be the compactification of $S \in \mathfrak{S}(k)$ with 
respect to the given metric $\rho$. A point $x \in \overline{S}$ 
is called a \emph{point at infinity} if $x \notin S$, and it is called
a \emph{point in $S$} or a \emph{finite point} otherwise.  A rectifiable curve $\gamma : I \to S$, 
where $I$ is one of the intervals $(0,1), [0,1), (0,1]$ or $[0,1]$, is called a \emph{locally shortest} curve 
if its extension $\overline{\gamma} : \overline{I} \to \overline{S}$ defines a shortest curve in 
$\overline{S}$. A locally shortest curve
with ends $x, y \in \overline{S}$ will be denoted by $[x,y]$ as in the case
of shortest curves.

\begin{definition}\label{D:gentri}
A \emph{generalized simple geodesic triangle}, or just a triangle for brevity,
on an Aleksandrov surface $S \in \mathfrak{S}(k)$ is a closed (not necessarily compact) subset $\D$ of $S$
such that its interior is homeomorphic to the open unit disc and it satisfies
the following properties:
\begin{enumerate}[(a)]
 \item the boundary of $\D$ consists of three locally shortest curves, called
       \emph{sides} of $\D$;
 \item each pair of sides of $\D$ has one and only one common endpoint in
       $\overline{S}$, called a \emph{vertex} of $\D$;
 \item the $\overline{S}$-closure $\overline{\D}$ of $\D$ is a domain $R_k$.
\end{enumerate}
\end{definition}

Roughly speaking, a triangle is a simple triangle (as defined in Section~\ref{S:Fact}) 
with some vertices deleted, and we require that it is also small to be a domain $R_k$. 
Also note that the condition (c) above makes it possible to define
an angle of $\D$ at a vertex at \emph{infinity} by the formula \eqref{upperangle}. 

Throughout this paper \emph{triangles in $S$} always mean the sets 
described in Definition~\ref{D:gentri}, 
while \emph{triangles in the ``model domain'' $S(k)$} are the
usual \emph{compact} triangles in the Euclidean plane, hemispheres or hyperbolic planes.

\subsection{Total angle and angular curvature}\label{theorem}
In this subsection we discuss the meaning of total angle and angular curvature.
To do this, let $S \in \mathfrak{S}(k)$ with a tiling $\T$. For
a vertex $v \in \overline{S}$ of a triangle $\D \in \T$, we denote by 
$\theta(v,\D)$ the angle of $\D$ at $v$ as before, and define the \emph{total angle} 
$T(v)$ of $v$ to be 
\begin{equation}\label{E:TA}
 T(v) := \sum_{\D \in \T} \theta(v, \D),
\end{equation}
with the convention that $\theta(v, \D) = 0$ if $v$ is not a vertex of $\D$.
By the Gauss-Bonnet formula (\cite{AZ}, p.~214) one can check that
this definition of total angle coincides with \eqref{ta} for a \emph{finite}
vertex $v$. Also note that \eqref{E:TA} is valid for a
vertex $v$ at infinity by defining $T(v) = \infty$. In fact if we assume the condition (R1) of Theorem~\ref{T}, 
then because $S$ is simply-connected
the case $T(v) < \infty$ cannot happen for a vertex $v$ at infinity, unless $\T$ contains only finitely 
many triangles. 
But then $S \cup \{v\}$ becomes compact, hence $S$ is parabolic.
In addition, one may check that $\sum_{\D \in \T} K(\D) = 4 \pi$ in this case and 
the condition (M2) of Theorem~\ref{T} is never satisfied. So without loss of generality, we will always assume 
that $T(v) = \infty$ for a vertex $v$ at infinity. 

The reason why we consider the angular curvature of a triangle is the following.
For a given triangle $\D \in \T$, we distribute to $\D$ the 
integral curvature concentrated at its 
vertices, weighted by the angle at each vertex. Thus if $\D$ has 
vertices $v_i$ and sides
$L_i$, $i=1,2,3$, the integral 
curvature of the \emph{closed} triangle $\D$ would be (with the interpretation
$\infty / \infty = 1$ and $a / \infty = 0$ for $a < \infty$)
\begin{equation}\label{E:Intcur2}
 \sum_{i=1}^3 \frac{2\pi - T(v_i)}{T(v_i)} \theta(v_i,\D) 
   + \omega(\D^\circ) + \sum_{i=1}^3 \tau(L_i),
\end{equation}
where $\tau(L_i)$ is the so-called \emph{left turn} (cf. \cite{AZ} or \cite{Res}) of the shortest curve
$L_i$. (The reason we have in \eqref{E:Intcur2} the left turn of sides of $\D$ is the same as before:
we want to compute the integral curvature of a \emph{closed} triangle.) 
On the other hand, the Gauss-Bonnet formula
applied to a triangle $\D$  implies that
\begin{equation}\label{E:GB}
 \sum_{i=1}^3 \theta(v_i,\D) - \pi - \sum_{i=1}^3 \tau(L_i) = \omega(\D^\circ).
\end{equation}
Now one can easily see 
that the quantity in \eqref{E:Intcur2} is same as the angular curvature 
$K(\D)$ defined in \eqref{ac}. 

\subsection{\"Uberlagerungsfl\"ache}\label{uber}
Let $X$ and $Y$ be two-dimensional topological manifolds and $\varphi : X \to Y$
a topologically holomorphic\footnote{continuous, open and discrete} map. 
Then by a theorem of Stoilov \cite{St}, for every $x \in X$ there are a positive integer $m$, called the local degree of $\varphi$ at $x$, 
and complex local coordinates $z$ and $w$ in neighborhoods of $x$ and $\varphi (x)$ such that $w(\varphi(x)) = z(x)^m$. 
Therefore if $Y$ is a Riemann surface, there exists a unique conformal structure on $X$ which makes $\varphi$ 
holomorphic, and
the pair $(X, \varphi)$ is called  a \emph{Riemann surface spread over} $Y$ (\"uberlagerungsfl\"ache in German). Two such pairs $(X_1, \varphi_1)$ and
$(X_2, \varphi_2)$ are called equivalent if there is a homeomorphism $h : X_1 \to X_2$ such that 
$\varphi_1 = \varphi_2 \circ h$, and strictly speaking, a Riemann surface spread
over $Y$ is an equivalence class of such pairs.

Now suppose $(X,\varphi)$ is an open simply-connected 
Riemann surface spread over $Y = \mathbb{C}$, $\overline{\mathbb{C}}$ or $\mathbb{D}$.
Then the surface $X$ equipped with the pull-back metric, i.e., the metric whose length element is of the form
\begin{equation}\label{metric}
 \frac{2 |\varphi'(z)||dz|}{1+\chi |\varphi(z)|^2}
\end{equation}
with $\chi = 0, 1$ or $-1$ depending on whether $Y$ is $\mathbb{C}$, $\overline{\mathbb{C}}$ or $\mathbb{D}$, respectively,
is a special type of an Aleksandrov surface, called a surface with \emph{polyhedral} metric \cite{Res}.
In this connection the study of Aleksandrov surfaces naturally arises in the field of the function theory. Also note that
the Riemann surface $X$ equipped with the metric defined by \eqref{metric} is in $\mathfrak{S}(\chi)$.

\subsection{Notations}\label{N:symbols}
We always use the letter $\rho$ to denote the intrinsic metric of a given 
Aleksandrov surface $S$, or the length element of this metric. The area of 
a region in $S$ will be denoted by $|\cdot|_{\rho}$, the length of a curve 
in $S$ by $\ell_{\rho}(\cdot)$, and the distance between two points $x$ and
$y$ in $S$ by $d_{\rho}(x,y)$. When we consider the metric $\lambda$, which will be
constructed in Section~\ref{S:Proof}, the subscript $\rho$ will 
be modified to $\lambda$.

\section{Three lemmas}\label{trilemmas}
We first state several theorems of A.~D.~Aleksandrov, which will
be used in the proof of Theorem~\ref{T}.

\begin{theorem}[Aleksandrov \cite{Alex2}]\label{TB}
 The upper angles $\overline{\alpha}, \overline{\beta}, \overline{\gamma}$ 
 of an arbitrary triangle $\D$ in a domain $R_k$ are
 not greater than the corresponding angles $\alpha_k, \beta_k, \gamma_k$ of
 the associated triangle $\D_k$ on $S(k)$.
\end{theorem}

Note that this theorem says that each angles of $\D$ are less than 
or equal to the corresponding angles of $\D_k$, while the second condition of
Definition~\ref{domainrk} says that only the sum of angles of $\D$ is at most
that of $\D_k$. 

\begin{theorem}[Aleksandrov \cite{Alex2}]\label{Alek1}
Any two points of $R_k$ are joined by a unique shortest curve.
\end{theorem}

\begin{theorem}[Aleksandrov \cite{Alex1}]\label{TC}
Let $P$ be a polygon in $S \in \mathfrak{S}(k)$, $k>0$. If the
perimeter of $P$ is less than $2 \pi k^{-1/2}$ and every two points in $P$ can
be joined by a unique shortest curve in $P$, then the area of $P$ does not exceed 
that of the polygon $P_0$ which has sides of the same length and is 
inscribed in a circle on a hemisphere of curvature $k$.
\end{theorem}

Note that  Theorems \ref{Alek1} and \ref{TC} 
in particular imply that any triangle $\D \subseteq S \in \mathfrak{S}(k)$, $k > 0$, 
has area at most $2 \pi k^{-1}$. 

In the following two sections, $S$ will represent an Aleksandrov surface of curvature at most $k$
satisfying all properties of Theorem~\ref{T}. All notations in Theorem~\ref{T} (such as
$\T$, $\PP$) and in Section~\ref{N:symbols} will be used without 
further remarks.

\begin{lemma}\label{L:length}
Let $\D$ be a triangle in $\T$ and $\Gamma$ a Jordan arc contained 
in $\D$ with endpoints on $\partial \D$. Among two subarcs of $\partial
\D$ (or $\partial \overline{\D}$ if $\D$ has some vertices at infinity) 
divided by $\Gamma$, we take the subarc $\Lambda$ which contains the
smaller number of vertices. (If both subarcs contain the same number
of vertices, we choose either.) Then there exists a constant $C=C(\epsilon)$
such that
\[
 \ell_\rho(\Lambda) \leq C \ell_\rho(\Gamma).
\]
\end{lemma}

\begin{proof}
Because every side of $\D$ is a locally shortest curve, there is nothing to prove if
the endpoints of $\Gamma$, $x$ and $y$, are on a same side of $\D$. So we 
assume that $x$ and $y$ are on different sides of $\D$. In this case, $\Lambda$
must contain one and only one vertex of $\D$, denoted by $z$, in its interior
arc.

Let $\wD$ be the triangle contained in $\D$ with 
vertices $x$, $y$ and $z$. We denote
by $a$, $b$ and $c$ the length of sides of $\wD$ which are opposite to $x$,
$y$ and $z$, respectively, and by $\alpha$, $\beta$ and $\gamma$ the 
angles of $\wD$ at $x$, $y$ and $z$, respectively. Now if we show that 
\begin{equation}\label{E:maxtri}
 \max \{ a, b \} \leq C(\epsilon) \cdot c
\end{equation}
for some constant $C(\epsilon) > 0$, then we have proved the lemma since
\begin{align*}
 \ell_\rho(\Lambda) & = d_\rho(x,z) + d_\rho(y,z) = a + b \leq 2C(\epsilon)
 \cdot c\\
  & = 2 C(\epsilon) d_\rho(x,y) \leq 2 C(\epsilon) \ell_\rho(\Gamma).
\end{align*}

To show \eqref{E:maxtri} we first note that it can be assumed that 
$k \in \{0,1\}$.
In fact, we may take $k=0$ if $k \leq 0$ because $S \in \mathfrak{S}(k) \subseteq \mathfrak{S}(0)$.
If $k>0$, then a surface $S \in \mathfrak{S}(k)$ with the given metric $\rho$ can
be realized as a surface in $\mathfrak{S}(1)$ with the metric $\sqrt{k}\rho$, so in this
case $k$ can be assumed to be 1. 

We next construct a new triangle $\wD_k$ on $S(k)$ whose side-lengths 
are those of $\wD$. Theorem~\ref{TB}
then yields that the corresponding angles 
$\alpha_k$, $\beta_k$ and $\gamma_k$ are greater than or equal to $\alpha$,
$\beta$ and $\gamma$, respectively. In particular, $\gamma_k \geq \gamma \geq \epsilon$
since $\gamma$ is an angle of $\D \in \T$ and $\T$ satisfies the condition (R1) of Theorem~\ref{T}. 
Furthermore, one may assume that $\gamma_k \leq \pi/2$ by dividing
$\wD_k$ into two triangles with the bisector of $\gamma_k$ if necessary, 
because our goal here is to show \eqref{E:maxtri} in the triangle $\wD_k$. Thus
\begin{equation}\label{E:sinegamma}
 \sin \gamma_k \geq \sin \epsilon.
\end{equation}

Now if $k=0$, the Law of Sines and \eqref{E:sinegamma} imply
\[
 b = \frac{\sin \beta_k}{\sin \gamma_k} c \leq \frac{1}{\sin \gamma_k} c
  \leq \frac{1}{\sin \epsilon} c,
\]
and similarly $a \leq (\sin \epsilon)^{-1}c$, as desired. If $k=1$, 
then since the perimeter of $\D$ is bounded by $2\pi - \epsilon$ 
by the assumption (R2) with $k=1$, we have $0 \leq b \leq \pi - \epsilon/2$, or
\begin{equation}\label{E:sineb}
 b \leq \frac{\pi - \epsilon/2}{\sin (\pi - \epsilon/2)} \sin b
   \leq \frac{\pi}{\sin(\epsilon/2)}\sin b.
\end{equation}
On the other hand, by the Spherical Law of Sines (cf. \cite{Chau}, 
p.~151)
\[
 \frac{\sin a}{\sin \alpha}
 = \frac{\sin b}{\sin \beta} = \frac{\sin c}{\sin \gamma}
\]
and \eqref{E:sinegamma} we have 
\begin{equation}\label{E:sinb2}
 \sin b = \frac{\sin \beta_k}{\sin \gamma_k}\sin c \leq \frac{1}{\sin \epsilon}
 \sin c \leq \frac{c}{\sin \epsilon}.
\end{equation}
Therefore, combining \eqref{E:sineb} and \eqref{E:sinb2},
$b \leq \pi (\sin (\epsilon/2)\sin (\epsilon))^{-1} c$, and similarly
$a \leq \pi (\sin (\epsilon/2)\sin (\epsilon))^{-1} c$. This completes the proof.
\end{proof}

The proof of this lemma shows that $\max \{ a,b \} \leq C(\epsilon) \cdot c$ if 
$\gamma \geq \epsilon$
and the perimeter of $\wD$ is at most $(2\pi - \epsilon)k^{-1/2}$ for $k>0$. 
Then since $\min \{ \alpha, \beta, \gamma \} \geq \epsilon$ when $\wD = \D$, 
we have the following corollary.

\begin{corollary}\label{Cor}
For every $\D \in \T$ there exists a constant $C=C(\epsilon)$ such that $M(\D) \leq C m(\D)$, where
\begin{equation}\label{MM}
 \begin{aligned}
 M(\D) &:= \max \{\mbox{\emph{side-lengths of }} \D \}, \\
 m(\D) &:= \min \{\mbox{\emph{side-lengths of }} \D \}. 
 \end{aligned}
\end{equation}
\end{corollary}

We next have an isoperimetric inequality for \emph{small} domains.

\begin{lemma}\label{L:Iso}
 Suppose $D$ is an open set which is contained in a single triangle $\D \in \T$. 
  Then
 \begin{equation}\label{E:Iso}
  |D|_{\rho} \leq \frac{1}{2\pi} \ell_{\rho}(\partial D)^2.
 \end{equation}
\end{lemma}

\begin{proof}
Without loss of generality, we may assume that all vertices of $\D$
are finite, because otherwise the argument below holds with the $\overline{S}$-closure 
$\overline{\D}$ in place of $\D$.

Now if all vertices of $\D$ are finite, by considering each component of 
$D$ separately and by adding to $D$ 
all complementary components of $D$ which are compactly contained in $\D$, we may assume 
that $D$ is simply-connected. Since
$S \in \mathfrak{S}(k)$, we have (\cite{BZ}, p.~12)
\begin{equation}\label{E:Iso2}
 \ell_{\rho}(\partial D)^2 - 4\pi |D|_{\rho} + k |D|_{\rho}^2 \geq 0.
\end{equation}
This directly implies \eqref{E:Iso} when $k \leq 0$.

If $k>0$, then by the statement following Theorem~\ref{TC} we have 
\[
 |D|_{\rho} \leq |\D|_{\rho} \leq \frac{2\pi}{k}.
\]
This inequality together with \eqref{E:Iso2} implies 
\begin{equation*}
 2\pi |D|_{\rho} \leq (4\pi-k|D|_{\rho})|D|_{\rho} \leq 
 \ell_{\rho}(\partial D)^2,
\end{equation*}
as desired. 
\end{proof}

The following lemma is the only place where conditions (M1) and (M2) 
of Theorem~\ref{T} are used, 
even though these conditions are our main assumptions. 

\begin{lemma}\label{L:Comb}
 Suppose that $S$ is an Aleksandrov surface with a tiling $\T$
 and a partition $\PP$ such that the conditions \emph{(M1)} and 
 \emph{(M2)} in Theorem~\ref{T} hold. Let $\T'$ be a
 finite subset of $\T$ and $D := (\bigcup_{\D \in \T'} \D)^{\circ}$. Define
 \begin{align*}
  e_0 &:= \mbox{\emph{the number of edges on }} \partial D, \\
  f   &:= \mbox{\emph{the number of triangles in }}\overline{D}.
  \end{align*}
 Then 
 \[
 f \leq \frac{6 M^2}{\epsilon} e_0,
 \]
where $M$ and $\epsilon$ are the constants in Theorem~\ref{T}.
\end{lemma}

\begin{proof}
We first assume that $D$ consists of clusters in $\PP$; i.e., $D$ is of the form
 \begin{equation}\label{st}
  D = (\bigcup_{i =1}^N C_{\alpha_i})^{\circ}, 
 \end{equation}
for some $C_{\alpha_i} \in \PP$, $i = 1, \ldots ,N$. 
Then because the interior of each cluster in $\PP$ is connected, all components of 
$D$ and all the components of $(S \backslash \overline{D})$ which consist of finitely many triangles 
are also unions of clusters. That is, they are also expressed 
in the same way as the formula \eqref{st}. 
Therefore to prove the lemma in this case, 
we may assume that $D$ is simply-connected by adding to $D$ all the 
complementary components of $D$ consisting of finitely many triangles, 
and considering each components of $D$ separately.
 
 Let $V'$ be the set of vertices lying in the interior of $D$
 and let $v' := |V'|$. (Note that $V'$ does not contain any vertices at
 \emph{infinity} because $f$ is finite.)
 Then $e_0 + v'$ is the number of all vertices in the $\overline{S}$-closure $\overline{D}$, 
 since for simply-connected $D$ the number of edges on $\partial D$
 is same as that of vertices on $\partial \overline{D}$.
 Let $e$ be the number of edges of triangles in
 $\overline{D}$. Euler's formula then gives
  \begin{equation}\label{E:Euler}
   e_0 + v' + f = e + 1.
  \end{equation}

 On the other hand, each triangle has three edges and each 
 edge corresponds to 
 two different triangles except those on $\partial D$, each of which 
 corresponds to only one triangle. Thus
 \[
 3f + e_0 = 2e.
 \]
 Combining this with \eqref{E:Euler}, we get 
  \begin{equation}\label{E:SimEuler}
   e_0 = f - 2v' + 2.
  \end{equation}
  
 Now since $V'$ contains no vertices at infinity, we have
 $T(v) < \infty$ for all $v \in V'$. Therefore we have
 \begin{align*}
  2 \pi v' - \pi f &= 2 \pi \left( \sum_{v \in V'} \frac{T(v)}{T(v)}
  \right) - \pi f 
  = 2 \pi \sum_{v \in V'} \frac{1}{T(v)} \left( \sum_{\D \in \T'} \theta(v , \Delta)
     \right) - \pi f \\
  &\leq \sum_{\Delta \in \T'} \left( 2 \pi \sum_{i=1}^{3} 
        \frac{\theta(v_i , \Delta)}{T(v_i)} \right) - \pi f 
   = \sum_{\Delta \in \T'} K(\Delta) \leq - \epsilon \pi N,
 \end{align*}
 or
 \begin{equation*}
  2v' \leq f - \epsilon N.
 \end{equation*}
 This inequality and \eqref{E:SimEuler} imply
 \begin{equation}\label{sharp}
  e_0 \geq f - (f - \epsilon N ) = \epsilon N \geq \frac{\epsilon}{M} f
 \end{equation}
 because $f$ is less than or equal to $MN$.

 If $D$ is not a union of clusters, let $\{ C_1, \ldots , C_s \}$ 
 be the largest subset of 
 $\PP$ such that $C_i \cap \partial D$ contains an edge for all $i =
 1, \ldots , s$. Now let $\mathfrak{D} :=(D \cup \bigcup_{i=1}^s C_i)^{\circ}$,
 and note that $\mathfrak{D}$ has the form \eqref{st} because 
 the interior of every cluster in $\PP$ is connected. Define
 \begin{align*}
   \tilde{e_0} &:= \mbox{the number of edges on } \partial \mathfrak{D}, \\ 
   \tilde{f}   &:= \mbox{the number of triangles in } \overline{\mathfrak{D}}.
  \end{align*}
 
 Because $\partial D$ has $e_0$ edges and each of which corresponds to at 
 most two clusters, $s$ must be less than or equal to
 $2 e_0$. On the other hand, every
 edge on $\partial \mathfrak{D}$ must be an edge in $C_i$, for some 
 $i = 1, \ldots , s$. Furthermore each $C_i$ has at most $3M$ edges since it
 contains at most $M$ triangles. Therefore,
 \[
  \tilde{e_0} \leq 3Ms \leq 6M e_0.
 \]
 Since $f \leq \tilde{f} \leq (M / \epsilon) \tilde{e_0}$ by \eqref{sharp}, we get
 \[
  f \leq \frac{6M^2}{\epsilon} e_0
 \]
 as desired.
\end{proof}

\section{Proof of Theorem~\ref{T}.}\label{S:Proof}
For any $\D \in \T$, let $e$ be a side of $\D$. We define
\begin{equation*}
 \lambda_0 (z) := \frac{1}{\ell_{\rho}(e)} \, , \qquad \mbox{for all }
 z \in e \backslash \{\mbox{endpoints of } e \}.
\end{equation*}
We define $\lambda_0$ for the other sides of $\D$ in the same way and extend
$\lambda_0$ to the interior of $\D$ as the bounded solution of the Dirichlet problem
\begin{align*}
 \varDelta u = 0  & \qquad \mbox{in } \D^{\circ}, \\
 u(z)        = \lambda_0 (z) & \qquad \mbox{for } z \in \partial 
                          \D \backslash \{\mbox{vertices of } \D \}.
\end{align*}
We either leave $\lambda_0$ undefined on vertices or define it to be zero.
Finally we define the metric $\lambda$ as the metric whose length element
is $\lambda_0 \rho$. This length element is also denoted by $\lambda$.

Note that for any triangle $\D \in \T$ and any side $e$ of $\D$,
\begin{equation}\label{E:Leq1}
 \ell_{\lambda} (e) = 1
\end{equation}
by the definition of $\lambda$. Furthermore by Phragm\'en-Lindel\"{o}f's maximum 
principle (\cite{Ahl}, p.~38),
\begin{equation*}
 \frac{1}{M(\D)}\rho(z) \leq \lambda(z) \leq \frac{1}{m(\D)}\rho(z),
 \qquad \mbox{for a.e. } z \in \D,
\end{equation*}
where $M(\D)$ and $m(\D)$ are defined in \eqref{MM}. It follows that for any 
Borel set $D \subseteq \D$ and any rectifiable curve $\Gamma \subseteq \D$,
\begin{align}
 & \frac{1}{M(\D)^2}|D|_{\rho} \leq |D|_{\lambda} \leq 
 \frac{1}{m(\D)^2}|D|_{\rho}, \label{E:Dbound} \\
 & \frac{1}{M(\D)}\ell_{\rho}(\Gamma) \leq \ell_{\lambda}(\Gamma)
 \leq \frac{1}{m(\D)}\ell_{\rho}(\Gamma). \label{E:Lbound}
\end{align}

\begin{lemma}\label{L:Liso}
Suppose $D$ is an open set in $S$ which is contained in a triangle
$\D \in \T$. Then $|D|_{\lambda} \leq C \ell_{\lambda} (\partial D)$
for some $C = C(\epsilon)$.
\end{lemma}

\begin{proof}
First we assume that $\ell_{\lambda} (\partial D) \geq 1$. Then by 
Lemma~\ref{L:Iso}, Corollary~\ref{Cor} and \eqref{E:Leq1}--\eqref{E:Lbound},
\begin{equation}\label{E:DleqL}
\begin{aligned}
 |D|_{\lambda} & \leq |\D|_{\lambda} \leq \frac{1}{m(\D)^2}|\D|_{\rho} \leq
      \frac{1}{m(\D)^2}\frac{1}{2\pi}\ell_{\rho}(\partial \D)^2 \\
 & \leq \frac{M(\D)^2}{2\pi m(\D)^2}\ell_{\lambda}(\partial \D)^2
      \leq \frac{9C^2}{2\pi} \leq \frac{9C^2}{2\pi}
      \ell_{\lambda} (\partial D), 
\end{aligned}
\end{equation} 
as desired. If $\ell_{\lambda} (\partial D) < 1$, a similar calculation
shows
\begin{align*}
 |D|_{\lambda} & \leq \frac{1}{m(\D)^2}|D|_{\rho} \leq
      \frac{1}{m(\D)^2}\frac{1}{2\pi}\ell_{\rho}(\partial D)^2 \\
 & \leq \frac{M(\D)^2}{2\pi m(\D)^2}\ell_{\lambda}(\partial D)^2
      \leq \frac{C^2}{2\pi}\ell_{\lambda}(\partial D), \notag
\end{align*}
and the lemma follows.
\end{proof}

In \eqref{E:DleqL}, we showed that for any $\D \in \T$, 
$|\D|_{\lambda} \leq (9C^2)/(2\pi)$.
Therefore by \eqref{E:Leq1} and Lemma~\ref{L:Comb},
we deduce the following corollary.

\begin{corollary}\label{Cor2}
As in Lemma~\ref{L:Comb}, let $D$ be an open set consisting of a finite number of triangles in $\T$.
Then there exists a constant $C=C(\epsilon, M)$ such that
\begin{equation*}
|D|_{\lambda} \leq C \ell_{\lambda} (\partial D).
\end{equation*}
\end{corollary}

Now suppose that a Jordan region $D$ in $S$ is given and $\D$ is a triangle in $\T$
such that $\partial D \cap \D^{\circ} \ne \emptyset$ and $D \nsubseteq \D$. Then
$\Gamma := \D^{\circ} \cap \partial D$ is a countable union of Jordan 
arcs in $\D$ with endpoints on $\partial \D$; i.e., $\Gamma = \bigcup_j 
\Gamma_j$ where $\Gamma_i \cap \Gamma_j = \emptyset$ if $i \ne j$, and each $\Gamma_j$
is a Jordan arc in $\D$ with endpoints on $\partial \D$.

For each $\Gamma_j$, let $\Lambda_j$ be one of the two closed subarcs of 
$\partial \overline{\D}$ which contains the smaller number of vertices, as in 
Lemma~\ref{L:length}, and $U_j$ the subregion of $\D$ enclosed 
by $\Gamma_j$ and $\Lambda_j$. Note that in our notation, $D$ and each $U_j$
are open regions, each $\Gamma_j$ is an open arc, while each $\Lambda_j$ and $\D$ 
are closed sets. Finally let's define $\alpha := \partial \D \cap D$ and
$\beta := \partial \D \backslash \overline{D}$.

\begin{lemma} \label{L:Either}
Either $\alpha \subseteq \bigcup_j \Lambda_j$ or 
$\beta \subseteq \bigcup_j \Lambda_j$.
If $\beta \nsubseteq \bigcup_j \Lambda_j$, then $D \cap
\Delta \subseteq \bigcup_j \overline{U}_j$.
\end{lemma}

\begin{proof}
Suppose that the first statement is not true. Then there exist 
$x$ and $y$ such that $x \in \alpha \backslash \bigcup_j \Lambda_j$ and 
$y \in \beta \backslash \bigcup_j \Lambda_j$. 
Since $x \in D$ and $y \in S \backslash \overline{D}$, 
we can find a Jordan arc $\Gamma_j$ which separates $x$ and $y$ in $\Delta$. Then either
$x \in \Lambda_j$ or $y \in \Lambda_j$ by the definition of $\Lambda_j$,
but this contradicts our assumption. 

Now suppose that $\beta \nsubseteq \bigcup_j \Lambda_j$. Then there exists
$y \in \beta \backslash \bigcup_j \Lambda_j$.
 So for any $x \in D \cap
\Delta$, either $x \in \overline{U}_j$ or $y \in \overline{U}_j$ for some $j$, by the
same argument as above. But $y \notin \Lambda_j 
= \overline{U}_j \cap \partial \Delta$, 
which implies that $x \in \overline{U}_j \subseteq \bigcup_j \overline{U}_j$. So
$D \cap \Delta \subseteq \bigcup_j \overline{U}_j$, as desired.
\end{proof}

If $\beta \nsubseteq \bigcup_j \Lambda_j$, we have $\partial (D \backslash \D)
\subseteq \partial D \cup \alpha$. In fact, if $x \in \partial (D \backslash
\D) \backslash \partial D$, then $x \in (\overline{D \backslash \D}) 
\subseteq \overline{D} \cap (\overline{S \backslash \D})$, $x \notin
D \backslash \D$ and $x \notin \overline{D} \backslash D$. In particular,
$x \in \overline{D}$, $x \notin \overline{D} \backslash D$ and $x \notin
D \backslash \D$. So $x \in D \cap \D$. But since $x \in
(\overline{S \backslash \D})$, we conclude that $x \in D \cap \partial \D
= \alpha$.
In this case, we claim that for some constant $C = C(\epsilon)$, 
\begin{align}
 & \ell_{\lambda}(\partial(D \backslash \D)) \leq
 \ell_{\lambda}(\partial D) + \ell_{\lambda}(\alpha) \leq
 \ell_{\lambda}(\partial D) + C \ell_{\lambda}(\Gamma), 
 \label{E:LL} \\
 & |D \backslash \D|_{\lambda} \geq |D|_{\lambda} - |D \cap \D|_{\lambda}
 \geq |D|_{\lambda} - C \ell_{\lambda}(\Gamma). \label{E:AA}
\end{align}
In fact, since $\beta \nsubseteq \bigcup_j \Lambda_j$ we have $\alpha \subseteq \bigcup_j \Lambda_j$ 
by Lemma~\ref{L:Either}.
Therefore by \eqref{E:Lbound} and Lemma~\ref{L:length}, we have
\begin{equation}\label{E:LCL}
\begin{aligned}
 \ell_{\lambda}(\alpha) &\leq
 \sum_j \ell_{\lambda}(\Lambda_j) \leq 
 \sum_j \frac{1}{m(\D)} \ell_{\rho}(\Lambda_j) 
 \leq C \sum_j \frac{1}{m(\D)} \ell_{\rho}(\Gamma_j) \\
 & \leq C \sum_j \frac{M(\D)}{m(\D)} \ell_{\lambda}(\Gamma_j) \leq
    C^2 \ell_{\lambda}(\Gamma), 
\end{aligned}
\end{equation}
and \eqref{E:LL} follows. Similarly by Lemma~\ref{L:Either},
Lemma~\ref{L:Liso} and the same estimate for 
$\sum_j \ell_{\lambda}(\Lambda_j)$,
\begin{align*}
 |D \cap \D|_{\lambda}  & \leq \sum_j |U_j|_{\lambda} 
  \leq \sum_j C \ell_{\lambda}(\partial U_j) \\ 
  & = C \sum_j [\ell_{\lambda}(\Lambda_j) + \ell_{\lambda}(\Gamma_j)]
  \leq (C^3 + C) \ell_{\lambda} (\Gamma),
\end{align*}
and \eqref{E:AA} is also proved.

Next let's consider the case $\beta \subseteq \bigcup_j \Lambda_j$. In this
case, we have $\partial (D \cup \D) \subseteq \partial D \cup \beta$. In fact,
if $x \in \partial (D \cup \D) \backslash \partial D$, then $x \in 
(\overline{D \cup \D}) = \overline{D} \cup \D$, $x \notin (D \cup \D)^{\circ} 
\supseteq D \cup \D^{\circ}$ and $x \notin \overline{D} \backslash D$. In
particular, $x \notin \overline{D} \backslash D$ and $x \notin D$, so
$x \notin \overline{D}$. Hence $x \in \D$. But $x \notin \D^{\circ}$, so we 
must have $x \in \partial \D$; i.e., $x \in \partial \D \backslash \overline{D}
= \beta$.

Then since $\beta \subseteq \bigcup_j \Lambda_j$, the estimate for
$\sum_j \ell_{\lambda} (\Lambda_j)$ in \eqref{E:LCL} implies that
\begin{equation}\label{E:LLL}
 \ell_{\lambda} (\partial (D \cup \D)) \leq \ell_{\lambda} (\partial D)
 + \ell_{\lambda} (\beta) \leq \ell_{\lambda} (\partial D) +
 C \ell_{\lambda} (\Gamma)
\end{equation}
for some $C = C(\epsilon)$. Also note that
\begin{equation}\label{E:AAA}
 |D \cup \D|_{\lambda} \geq |D|_{\lambda} \geq 
 |D|_{\lambda} - C \ell_{\lambda} (\Gamma).
\end{equation}
Collecting \eqref{E:LL}, \eqref{E:AA}, \eqref{E:LLL} and \eqref{E:AAA},
we have the following lemma.

\begin{lemma} \label{L:Final}
If $D$ is a Jordan region in $S$ and $\D$ a triangle in $\T$ such that
$\D^{\circ} \cap \partial D \ne \emptyset$ and $D \nsubseteq \D$, then we may obtain
a region $D'$ by properly
adding or subtracting $\D$ from $D$ so that 
the following properties hold for some constant $C=C(\epsilon)$:
\begin{align}
 & \ell_{\lambda}(\partial D') \leq \ell_{\lambda} (\partial D) + 
  C \ell_{\lambda} (\D^{\circ} \cap \partial D), \label{E:FL} \\
 & |D|_{\lambda} \leq |D'|_{\lambda} + C \ell_{\lambda} (\D^{\circ} \cap 
 \partial D). \label{E:FA}
\end{align}
\end{lemma}

Note that $D$ does not have to be a Jordan region in Lemma~\ref{L:Final}.
The conclusion of this lemma is still valid under a weaker assumption that 
$D$ is open and $\D^{\circ} \cap \partial D$ is a countable disjoint union of Jordan 
arcs in $\D$ with endpoints on $\partial \D$.

The following Theorem is due to L.~Ahlfors.

\begin{theorem}[Ahlfors]\label{TD} 
Let $S$ be an open simply-connected Riemann surface with a conformal metric
$\lambda$. If $\lambda$ allows a linear isoperimetric inequality for every 
Jordan region in $S$, then $S$ is hyperbolic; i.e., conformally equivalent to 
the unit disc.
\end{theorem}

For the proof of Theorem~\ref{TD}, see for example \cite{Hay},
pp.~143--144. In fact, this theorem is proved in \cite{Hay} only when 
$\lambda$ is the induced spherical metric, but one can check that the same 
argument works when $\lambda$ is merely \emph{locally integrable}.

Note that the surface $S$, when it is equipped with the metric $\lambda$ constructed in this section,
does not have to be an Aleksandrov surface. However, 
one can see that Theorem~\ref{TD} is still applicable to this case,
because by \eqref{hu} the pair $(S, \rho)$, where $\rho$
is the given metric which makes $S$ an Aleksandrov surface,
can be identified with a Riemann surface  and a conformal metric of the form \eqref{eqre},
and  $\lambda$ is conformally equivalent to $\rho$ by the construction.

\begin{proof}[Proof of Theorem~\ref{T}] 
By Theorem~\ref{TD}, it suffices to
show a linear isoperimetric inequality for every Jordan region $D$ in $S$.
If $D$ is contained in a single triangle in $\T$ or is a union of triangles, then
there is nothing to prove because of Lemma~\ref{L:Liso} or Corollary~\ref{Cor2},
respectively. So we assume that there exist a finite number of triangles
$\D_1, \ldots, \D_N$ in $\T$ such that $\D_j^{\circ} \cap \partial D \ne \emptyset$ and $D \nsubseteq \D_j$
for $j = 1, \ldots, N$. Since $\partial D$ is compact, there are only
finitely many such triangles.

Let $D_0 := D$. Then we successively construct $D_j$ by either adding 
or subtracting $\D_j$ from $D_{j-1}$ so that \eqref{E:FL} and 
\eqref{E:FA} are satisfied. (The $D_j$'s are not Jordan regions, but
since $\D_j^{\circ}$'s are disjoint we can still apply Lemma~\ref{L:Final}; 
see the statement following Lemma~\ref{L:Final}.) By the construction, $D_N$ 
is a union of triangles in $\T$. Moreover we have 
$\D_j^{\circ} \cap \partial D_{j-1} =  \D_j^{\circ} \cap \partial D$,
because the $\D_j^{\circ}$'s are disjoint. Therefore by Lemma~\ref{L:Final} and Corollary~\ref{Cor2}, 
\begin{align*}
|D|_{\lambda} &= |D_0|_{\lambda} \leq |D_1|_{\lambda} + 
    C \ell_{\lambda} (\D_1^{\circ} \cap \partial D) \\
&\leq |D_2|_{\lambda} +
    C \ell_{\lambda} (\D_2^{\circ} \cap \partial D) +
    C \ell_{\lambda} (\D_1^{\circ} \cap \partial D) \leq \cdots \\
&\leq |D_N| + 
    C \sum_{j=1}^N \ell_{\lambda} (\D_j^{\circ} \cap \partial D) \leq
    C \ell_{\lambda} (\partial D_N) + C \ell_{\lambda} (\partial D) \\
&\leq C [\ell_{\lambda} (\partial D_{N-1}) + 
    C \ell_{\lambda} (\D_N^{\circ} \cap \partial D)] +
    C \ell_{\lambda} (\partial D) \\
&\leq C [\ell_{\lambda} (\partial D_{N-2}) + 
    C \ell_{\lambda} (\D_{N-1}^{\circ} \cap \partial D) +
    C \ell_{\lambda} (\D_N^{\circ} \cap \partial D)] +
    C \ell_{\lambda} (\partial D) \\
&\leq \cdots \leq C \ell_{\lambda} (\partial D_0) +
    C^2 \sum_{j=1}^N \ell_{\lambda} (\D_j^{\circ} \cap \partial D) +
    C \ell_{\lambda} (\partial D) \\
& \leq (C^2 +2C) \ell_{\lambda} (\partial D),
\end{align*}
where $C=C(\epsilon, M)$. This completes the proof.
\end{proof}
   
\begin{remark}\label{rem}
Note that in the proof of Theorem~\ref{T} the assumption $S \in \mathfrak{S}(k)$ is never used 
globally---what we only used is the fact that every triangle $\D$ (or $\overline{\D}$) is a domain $R_k$.
Also note that all constants in the proof of Theorem~\ref{T} are independent of $k$.
(If we replace the assumption $(2\pi - \epsilon)k^{-1/2}$ in (R2) by
$2\pi k^{-1/2} - \epsilon$, then the constant in Lemma~\ref{L:length} 
will depend on $k$.) 
From these observations, we deduce the following stronger result:
suppose that an Aleksandrov 
surface $S$ has a tiling $\T$ and a partition $\PP$ which satisfies 
all assumptions in Theorem~\ref{T} with $k=k(\D) < \infty$. Then $S$ is hyperbolic. 
\end{remark}

\begin{remark}\label{rkk}
When $k>0$, suppose that we have
\medskip

(R1)$'$ (the length of every side of $\D$) $\geq \epsilon$, for every $\D \in \T$

\medskip
\noindent instead of (R1) in Theorem~\ref{T}. Then one can see that our proof 
yields a linear isoperimetric inequality
with respect to the original metric $\rho$ since $|\D|_\rho \leq 2 \pi k ^{-1}$ and every side of
a triangle $\D \in \T$ has length at most $\pi k^{-1/2}$. (Also one may check that 
Lemma~\ref{L:length} is still valid in this case. In fact, (R1) is a result of (R1)$'$ and (R2) when 
the curvature is constant $k>0$, which is exactly what we needed for the proof of Lemma~\ref{L:length}.)
When $k \leq 0$, the same is true using

\medskip
(R1)$''$ $\epsilon \leq$ (the length of every side of $\D$) $\leq M$, for every $\D \in \T$.

\medskip
\noindent Therefore in these cases the Aleksandrov surface satisfies a linear isoperimetric inequality
with respect to $\rho$, hence it is also \emph{Gromov hyperbolic} \cite{CDP}, \cite{BE3}, \cite{Oh}.
\end{remark}

\section{Construction of Example~\ref{Ex}}\label{S:Ex}
For given two Riemann surfaces $(X_i, \pi_i)$, $i=1,2$, spread over the complex plane, 
we use the following notation $(X_1, \pi_1) \subseteq (X_2,\pi_2)$ if
$X_1 \subseteq X_2$ and $\pi_2 |_{X_1} = \pi_1$. A point $x \in X_1$ (or a set $A \subseteq X_1$) 
is called a point (or a set) \emph{over} $y \in \mathbb{C}$ (or $B \subseteq \mathbb{C}$)
if $\pi_1 (x) = y$ (or $\pi_1 (A) = B$, respectively). 

In this section, we construct a sequence 
$(S_1, \varphi_1) \subseteq (S_2, \varphi_2) \subseteq \cdots$
of Riemann surfaces (with boundary) spread over proper subsets of $\mathbb{C}$,
and the limit surface $(S, \varphi)$, $S= \bigcup_{n=1}^\infty S_n$ and $\varphi = \lim_{n \to \infty} \varphi_n$,
will serve as the surface described in Example~\ref{Ex}. In fact, the surface $S$ equipped 
with the pull-back Euclidean metric, i.e., the metric whose length
element is of the form $|\varphi' (z)||dz|$, will be the Aleksandrov surface satisfying all properties
in Example~\ref{Ex}. 

Before constructing this surface, however, we fix some more notations:
we will always take the convention that each $S_n$ as well as $S$ is equipped with the pull-back Eulidean metric, 
and terminology such as shortest curves, area, distance, etc., 
will be used in this context. The symbol $|\cdot|$ will denote the area of a region in $S_n$ with respect 
to this metric, and $d(\cdot  ,\cdot)$ the distance between two points in $S_n$. 
Next, the term \emph{anticlockwise} will be used for the preferred orientation
of $\partial S_n$, where the orientation of $S_n$ is defined so that $\varphi_n$ is an orientation preserving map.
If $(\Sigma_n, \tilde{\varphi}_n)$ is another Riemann surface spread over the plane such that
$(S_n, \varphi_n) \subseteq (\Sigma_n, \tilde{\varphi}_n)$, then the term \emph{anticlockwise} will
also be used for $\partial \Sigma_n$ in the same way. Finally, for given
$R > 0$, $B(R)$ is the ball of radius $R$ in $\mathbb{C}$ with center
at the origin, $\overline{B}(R)$ its closure, and $C(R)$ its boundary.

Now let us construct the example. We first take an equilateral triangle
$\D_0$ inscribed in $C(R_1)$, for a sufficiently small $R_1 > 0$  such that
$|\D_0| \leq \epsilon$. Let $S_1 := \D_0$, $\varphi_1 := identity$ and 
$\T_1 := \{\D_0\}$.

Suppose that we have constructed $(S_1, \varphi_1) \subseteq
(S_2, \varphi_2) \subseteq \cdots \subseteq (S_n, \varphi_n)$, their tilings
$\T_1 \subseteq \T_2 \subseteq \cdots \subseteq 
\T_n$, and positive real numbers $R_1 \leq R_2 \leq \cdots \leq R_n$ 
which satisfy the following properties:
\begin{enumerate}[(a)]
 \item each $(S_i, \varphi_i)$ is a surface spread over a proper subset of 
       $\overline{B}(R_i)$ for $i=1,2,\ldots, n$; 
 \item each $\T_i$ is a tiling of $S_i$, $i=1,2,\ldots,n$, consisting of finitely many triangles;
 \item each triangle in $\T_n$ has area at most $\epsilon$;
 \item all vertices of triangles in $\T_n$ lie over $\bigcup_{i=1}^n C(R_i)$; 
 \item a vertex $v$ belongs to $\partial S_i \subseteq S_n$ if and only if 
       $v$ is over $C(R_i)$, $i=1,2,\ldots,n$;
 \item the set of critical points of $\varphi_n$ is equal to the set of all
       the vertices of triangles in $\T_{n-1}$, and at each critical point $\varphi_n$ has the local degree 2.
\end{enumerate}

Let $t_i$ be the number of vertices of triangles in $\T_n$ which are over $C(R_i)$, $i=1,2,\ldots,n$, and 
choose $R_{n+1}$ such that
\begin{equation}\label{E:Rn}
 R_{n+1} \geq R_n \exp\{2\pi n \sum_{i=1}^n t_i +2 \pi n \}.
\end{equation} 
We then extend $(S_n, \varphi_n)$ to a branched covering $(\Sigma_n, \tilde{\varphi}_n)$ of 
$\overline{B}(R_{n+1})$ so that $(S_n, \varphi_n) \subseteq (\Sigma_n, \tilde{\varphi}_n)$
and the set of critical points of $\tilde{\varphi}_n$ is
equal to that of $\varphi_n$. For example, 
$\Sigma_1 = \overline{B}(R_2)$ and $\tilde{\varphi}_1 = identity$ since $\varphi_1$ has
no singular point, and $\tilde{\varphi}_2 : \Sigma_2 \to \overline{B}(R_3)$ is a four-sheeted branched covering  
with critical points of the local degree 2 at the vertices of $\D_0$.
Similarly one can check by the Riemann-Hurwitz formula that  
$\tilde{\varphi}_n : \Sigma_n \to \overline{B}(R_{n+1})$ is a ($\sum_{i=1}^{n-1} t_i +1$)-sheeted branched covering.
Now let $\{v_1, v_2, \ldots, v_{t_n} \}$ be an anticlockwise enumeration (modulo $t_n$) 
of the vertices of triangles in $\T_n$ which belong to
$\partial S_n$. Then choose $t_n$ points
$w_1, w_2, \ldots , w_{t_n}$ on $\partial \Sigma_n$, also enumerated 
anticlockwise in modulo $t_n$, so that $|\D^j_n| \leq \epsilon$, where $\D_n^j$ is the triangle
in $\Sigma_n$ with vertices $v_j, v_{j+1}$ and $w_{j+1}$. This is possible by making the angle of 
$\D_n^j$ at $v_{j+1}$ appropriately small. (See Figure 1).

\begin{picture}(190,160)(-60, -10)
 \put(0,10){\arc[30](100,10){80}}
 \put(0,10){\arc[50](180,10){40}}
 \drawline[100](31,105.8)(167,77)
 \drawline[100](51,96)(167,77)
 \drawline[100](51,96)(177,42.3)
 \drawline[100](68.3,83)(177,42.3)
 \drawline[100](31,105.8)(51,96)
 \drawline[100](51,96)(68.3,83)
 \put(22,113){$v_{j+1}$}
 \put(43,88){$v_j$}
 \put(67,73){$v_{j-1}$}
 \put(170,77){$w_{j+1}$}
 \put(182,40.3){$w_j$}
 \put(90,0){\textbf{Figure 1}}
\end{picture}

Now for each $j =1,2,\ldots, t_n$, let $\Gamma_j$ be a copy of
$\overline{B}(R_{n+1})$ and $\pi_j : \Gamma_j \to \overline{B}(R_{n+1})$ 
the identity map. On each $\Gamma_j$, we mark two points
$v^j$ and $w^j$ so that $\pi_j (v^j) = \tilde{\varphi}_n(v_j)$ and 
$\pi_j (w^j) = \tilde{\varphi}_n(w_j)$. 
Then for each $j$, we cut $\Gamma_j$ along the shortest curve $[v^j, w^j]$, 
and cut $\Sigma_n$ along the 
shortest curves $[v_j, w_j]$
for all $j=1,2,\ldots,t_n$. Now we glue each $\Gamma_j$ to $\Sigma_n$, 
$j=1,2,\ldots,t_n$, along the corresponding cuts. 
This is possible since
$\pi_j ([v^j, w^j]) = \tilde{\varphi}_n ([v_j, w_j])$ for all $j$, and in this way we have
a new Riemann surface $(\Sigma_n', \tilde{\varphi}_n')$ spread over $\overline{B}(R_{n+1})$. 

We transfer notations such as $w_j, v_j$, etc., to $\Sigma_n'$ in the 
following way: each of $w_j \in \partial \Sigma_n$ and $w^j \in \partial 
\Gamma_j$ corresponds to two points in $\partial \Sigma_n'$, which are also 
denoted by $w_j$ and $w^j$; we place them in the order $\ldots, 
w_j, w^j, w_{j+1}, \ldots$ anticlockwise; each of $v_j \in \Sigma_n$ and
$v^j \in \Gamma_j$ corresponds to only one
point in $\Sigma_n'$, which is denoted by $v_j$; $\D_n^j$ is regarded as a 
triangle in $\Sigma_n'$ as well as in $\Sigma_n$ with vertices $v_j, v_{j+1}$
and $w_{j+1}$.

To complete the construction of $(S_{n+1}, \varphi_{n+1})$, let $M_n$ be a
positive integer which will
be determined later, and choose a set of points 
\[
 \{w_j^k : 1 \leq j \leq t_n, 0 \leq k \leq M_n \},
\]
on $\partial \Sigma_n'$, where $w_j^k$'s are enumerated anticlockwise for fixed $j$, 
so that for all $j$ and $k$, $w_j^0 = w_j = w_{j-1}^{M_n}$ and 
$d(w_j^k, w_j^{k+1}) = d(w_j^{k-1}, w_j^k)$. 
Let $\D_n^{jk}$ be the triangle in $\Sigma_n'$ with
vertices $w_j^k$, $w_j^{k+1}$ and $v_{j}$ for $1 \leq j \leq t_n$ and
$0 \leq k \leq M_n-1$. Taking $M_n$ sufficiently large, we have
$|\D_n^{jk}| \leq \epsilon$. Then the set
\[
 \T_{n+1} := \T_n \cup \{\D_n^j : 1 \leq j \leq t_n\} \cup
      \{\D_n^{jk} : 1 \leq j \leq t_n, 0 \leq k \leq M_n-1 \}
\]
is a tiling of $S_{n+1} := \bigcup_{\D \in \T_{n+1}} \D$, 
and we define 
\[
\varphi_{n+1} = \tilde{\varphi}_n'|_{S_{n+1}}.
\] 
Note that $S_{n+1}$ is obtained from $\Sigma_n'$ by cutting off some 
parts outside triangles in
$\T_{n+1}$, and $(S_{n+1}, \varphi_{n+1})$, $\T_{n+1}$ and $R_{n+1}$ satisfy all 
properties (a)--(f) above.

Now let $S:= \bigcup_{n=1}^{\infty} S_n = \bigcup_{n=1}^{\infty} \Sigma_n'$, 
$\varphi = \lim_{n \to \infty} \varphi_n$ and $\T := \bigcup_{n=1}^{\infty} \T_n$. 
Clearly the pair $(S, \varphi)$ is a Riemann surface spread over the plane
with the tiling $\T$. Since each $S_n$ is simply 
connected, so is $S$. Furthermore, the area of all triangles in $\T$ is at 
most $\epsilon$ by (c) and each vertex has total angle $4\pi$ by (f). 
Hence it remains to show that $S$ is conformally equivalent to $\mathbb{C}$, 
i.e., parabolic.

Let $B_n := B(R_{n+1}) \backslash \overline{B}(R_n) \subseteq \mathbb{C}$ and  
$A_n := \varphi^{-1}(B_n) \cap \Sigma_n'$.
Then one can easily see that $\{A_n\}$ is a sequence of concentric 
annuli in $S$.
Moreover, since each $\Sigma_n'$ is a $(\sum_{i=1}^n t_i +1)$-sheeted branched
covering of $\overline{B}(R_{n+1})$ and  
there is no critical point of $\varphi$ over $B_n$,
each $A_n$ is a $(\sum_{i=1}^n t_i +1)$-sheeted unbranched covering of
$B_n$.
Then because the module (cf. \cite{Ahl} or \cite{Va}) of $B_n$ is
$(2\pi)^{-1} \log (R_{n+1}/R_n)$, the module of $A_n$ is
\[
\frac{1}{2\pi(\sum_{i=1}^n t_i + 1)} \log \frac{R_{n+1}}{R_n},
\]
which is greater than or equal to $n$ by \eqref{E:Rn}. 
So the module of $A_n$ tends to
$\infty$ as $n \to \infty$, and this shows that $S$ is parabolic. In fact,
a comparison theorem (\cite{Ahl}, p.~54) implies that the module of each
$A_n$ is less than or equal to that of the unbounded component of 
$S \backslash A_1$, which should be finite if $S$ were hyperbolic.  

\section{Line Complexes and Hyperbolicity}\label{S:Appl}
Suppose that $(S, \varphi)$ is a simply-connected Riemann surface spread over the sphere
such that for some $q \geq 2$ points $a_1, \ldots, a_q \in \overline{\mathbb{C}}$, the restriction map
\begin{equation}
 \varphi: S \backslash \{ \varphi^{-1}(a_j) : j=1,2,\ldots,q \} \to \overline{\mathbb{C}} \backslash \{a_1,\ldots,a_q\}
\end{equation}
is a topological covering map. This pair $(S, \varphi)$ is called
a \emph{Riemann surface of class $F_q$} ramified over the points
$a_1, \ldots, a_q$, and denoted by $(S, \varphi) \in F_q(a_1, \ldots, a_q)$, or $S \in F_q$ if there is no confusion. 

Now we draw through the points $a_\nu$
in the order $a_1, \ldots, a_q, a_1$ a Jordan curve $\gamma$, called a \emph{base curve},
decomposing $\overline{\mathbb{C}}$ into two simply-connected open regions $G_1$ and $G_2$,
where $G_1$ is on the left of $\gamma$.  
These ``half sheets'' $G_1$ and $G_2$ are called polygons for short, $a_1, 
\ldots, a_q$ vertices, and the subarcs of $\gamma$, $(a_1 a_2), (a_2 a_3), \ldots,
(a_q a_1)$, sides of the polygons.
Next we choose two points $\circ \in G_1$ and $\times \in G_2$, and for each $j =1,2, \ldots, q$, let $\gamma_j$ be a
simple arc joining these two points and crossing the side $(a_j a_{j+1})$
at exactly one point so that $\gamma_j \cap \gamma_k = \{\circ, \times \}$
for $i \neq k$.
Let $\Gamma_0$ be the graph in $\overline{\mathbb{C}}$ consisting
of $\circ$, $\times$ and $\gamma_j$'s.
Then the pull-back graph $\Gamma := \varphi^{-1}(\Gamma_0)$ is a planar graph which is properly embedded in $\mathbb{C}$, and its homeomorphic
equivalence\footnote{Two properly embedded planar graphs $\Gamma_1$ and $\Gamma_{2}$ are called homeomorphically 
equivalent if there is a homeomorphism $h: \mathbb{C} \to \mathbb{C}$ such that $h(\Gamma_1)=\Gamma_2$.} class is 
called \emph{the line complex} or \emph{the Speiser graph}
of degree $q$ corresponding to $S \in F_q$ (\cite{Ne}, Chap.~XI). 
The components of $\mathbb{C} \backslash \Gamma$ will be called faces of $\Gamma$,
and $\Gamma$ has the following properties:
\begin{enumerate}[(i)]
 \item $\Gamma$ is connected;
 \item $\Gamma$ is bipartite: i.e., all vertices are split into two disjoint
       subsets, say $V_\circ := \varphi^{-1}(\circ)$ and $V_\times := \varphi^{-1}(\times)$, and every edge joins a point
       in $V_\circ$ and a point in $V_\times$;
 \item every vertex has the same degree $q$: i.e., each vertex corresponds to $q$ different
       edges;
 \item edges can be labeled by $1,2, \ldots, q$ so that they are placed
       counterclockwise around the points in $V_\circ$ and clockwise around 
       the points in $V_\times$.
\end{enumerate} 

Conversely, suppose we are given a line complex $\Gamma$ of degree $q$: i.e., suppose $\Gamma$ is 
(a homeomorphic equivalence class of) a planar graph which is properly 
embedded in $\mathbb{C}$ and satisfies (i)--(iv). Then for a fixed
base curve $\gamma$ passing through $a_1, \ldots, a_q$ in this order, there
exists a unique $S \in F_q$ whose corresponding line complex is $\Gamma$. 
Therefore for fixed base points $a_1, \ldots, a_q$ and a base curve $\gamma$, there is a one-to-one correspondence between 
Riemann surfaces of class $F_q$ and line complexes of degree $q$.

Each vertex $p$ of $\Gamma$ corresponds to a half sheet---a connected 
component of $\varphi^{-1}(G_1)$ or $\varphi^{-1}(G_2)$---of $S$, and each
face of $\Gamma$ corresponds to a singularity of $\varphi$. In fact, each face 
is a $2m$-gon for some $m \in \mathbb{N} \cup \{ \infty \}$, 
and if $m < \infty$ the face of $\Gamma$ corresponds to a point in $S$ at which $\varphi$ 
has the local degree $m$, 
and if $m = \infty$ the face corresponds to a logarithmic singularity of $\varphi$.

The following definition is due to R.~Nevanlinna.

\begin{definition}\label{D:RamExc}
Let $p$ be a vertex of $\Gamma$ and let $f_1, \ldots, f_q$ be the faces of
$\Gamma$ with $p$ on their boundaries. If each face $f_i$ is $2m_i$-gon for
$i = 1, \ldots, q$, \emph{the excess} $E_p$ of $\Gamma$ at $p$ is defined as
\begin{equation}\label{E:Exc}
 E_p := \sum_{i=1}^q \frac{1}{m_i} -q +2.
\end{equation}
We interpret $1/m_i = 0$ if $m_i = \infty$ for some $i$.
\end{definition}

A Riemann surface $S$ of class $F_q$  
is called \emph{regularly ramified} if there is a real number $E$ such that $E =  E_p$
for every vertex $p$ of the corresponding line complex $\Gamma$.

\begin{theorem}[Nevanlinna \cite{Ne}]\label{T:Ne}
Suppose that $S$ is a regularly ramified open Riemann surface of class $F_q$. Then $S$ is 
parabolic if $E=0$ and  hyperbolic if $E<0$.
\end{theorem}

Our purpose in this section is to generalize the hyperbolic case $E<0$ of
this theorem. Also note that the case $E>0$ happens only when $S$ is compact; i.e., $S = \overline{\mathbb{C}}$. 

Suppose that $\Gamma_\alpha$ is a connected subgraph of $\Gamma$ and 
$V_\alpha \subseteq V_\circ \cup V_\times$
the vertex set of $\Gamma_\alpha$. We then 
identify $\Gamma_\alpha$ with $V_\alpha$, and terminology in the set 
theory such as union, intersection, disjoint, etc., will be used for subgraphs in this 
context. A set $\PP(\Gamma)$ of subgraphs is called \emph{a partition} of 
$\Gamma$ if each element $\Gamma_\alpha \in \PP(\Gamma)$
is connected as a subgraph of $\Gamma$, elements in $\PP(\Gamma)$ are
disjoint, and $\bigcup_{\Gamma_\alpha \in \PP(\Gamma)} \Gamma_\alpha
= \Gamma$. Let $\#(\Gamma_\alpha)$ be the number of vertices in 
$\Gamma_\alpha$.

\begin{theorem}\label{T2}
Suppose that $\Gamma$ is a line complex with the associated Riemann surface
$S$ of class $F_q$. If $\Gamma$ has a partition $\PP(\Gamma)$ such that for some constants
$\epsilon > 0$ and $M >0$,
\begin{align*}
 &\mbox{\emph{(M1)$'$ }} \# (\Gamma_\alpha) \leq M \mbox{ for every }
     \Gamma_{\alpha} \in \PP(\Gamma), \\
 &\mbox{\emph{(M2)$'$ }} \sum_{p \in \Gamma_{\alpha}} E_p 
     \leq -\epsilon \mbox{ for every }\Gamma_{\alpha} \in \PP(\Gamma), 
\end{align*} 
then $S$ is hyperbolic.
\end{theorem}

\begin{proof}
Because the surface $S$ equipped with the pull-back spherical metric 
is in $\mathfrak{S}(1)$ (cf. Section~\ref{uber}), it suffices to
construct a tiling $\T$ and a partition $\PP$ which satisfy the four
conditions in Theorem~\ref{T}.

Without loss of generality, we assume that the base curve $\gamma$ is
a union of a finite number ($\geq q$) of geodesic line segments in 
$\overline{\mathbb{C}}$. Furthermore, we may assume that the angle of $\gamma$
at each $a_j$ is $\pi$ for all $j=1,\ldots,q$, and that the two components 
$G_1$ and $G_2$ of
$\overline{\mathbb{C}} \backslash \gamma$ have the same areas $2 \pi$. We then construct 
a finite tiling $\T'$ of $\overline{\mathbb{C}}$ so that the
vertex set of triangles in $\T'$ contains all $a_j$'s and all the endpoints of geodesic 
line segments in $\gamma$, and that both $G_1$ and $G_2$ are unions of 
a finite number of triangles in $\T'$.

The tiling $\T'$ of $\overline{\mathbb{C}}$ 
induces a tiling $\T$ of $S$ via the map $\varphi : S \to \overline{\mathbb{C}}$ as follows:
for each $\D' \in \T'$, the closure of each component of $\varphi^{-1} (\D^\circ)$ is
a triangle $\D \in S$ satisfying the conditions in Definition~\ref{D:gentri}, hence we let $\T$ be the collection
of such triangles. One can easily check that $\T$ satisfies the definition of a tiling given in the introduction, 
and note that the $\overline{S}$-closure $\overline{\D}$ of $\D$ is isometric to $\D'$.
Similarly, the partition $\PP(\Gamma)$ of $\Gamma$ induces a partition $\PP$ of $\T$. In fact,
each vertex $p$ of $\Gamma$ corresponds to a half sheet---a connected 
component of $\varphi^{-1}(G_1)$ or $\varphi^{-1}(G_2)$. Hence each subgraph
$\Gamma_\alpha \in \PP(\Gamma)$ corresponds to a union of half sheets, 
say $C_\alpha$, which is connected since each $\Gamma_\alpha$ is connected. 
Finally, since both $G_1$ and $G_2$ consist of triangles in $\T'$, each $C_\alpha$
consists of triangles in $\T$. Therefore $C_\alpha$ is a cluster. Now one can easily check that all the other 
assumptions about $\PP$ follow from the assumptions about $\PP(\Gamma)$.

We verify that the tiling $\T$ and the partition $\PP$ defined above satisfy all the conditions
of Theorem~\ref{T}. First, the conditions (R1) and (R2) are trivially satisfied because for every $\D \in \T$,
the $\overline{S}$-closure $\overline{\D}$ of $\D$ is isometric to a triangle $\D'$ in $\T'$, and $\T'$ contains 
only finitely many triangles. Furthermore (M1)$'$ implies (M1) since there are only finitely
many triangles in both $G_1$ and $G_2$.

It remains to show that the condition (M2) holds. To do this,
fix a vertex $p$ of $\Gamma$, and let $H$ be the corresponding half sheet
of $S$,
$f_1,\ldots,f_q$ the faces of $\Gamma$ with $p$ on their boundary, 
and $\nu_1, \ldots, \nu_q$ the corresponding singularities of $S$ lying on 
the $\overline{S}$-boundary $\partial H$ of $H$.
If $f_j$ is $2 m_j$-gon for $j=1,2, \ldots , q$, 
we have $T(\nu_j) = 2 \pi m_j$. Moreover, since the angle of $\gamma$ at 
$a_j$ is $\pi$, the angle of $\partial H$ at $\nu_j$ is also $\pi$. Finally, 
we have $T(v) = 2 \pi$ for every vertex $v \neq \nu_j$, $j=1,2, \ldots , q $, of triangles in $H$.
Therefore,
\begin{equation}\label{E:SumCur}
 \begin{aligned}
  \sum_{\D \subseteq H} & \sum_{i=1}^3 \frac{2\pi - T(v_i(\D))}{T(v_i(\D))}
   \theta(v_i(\D), \D) 
  = \sum_{j=1}^q  \sum_{\D \subseteq H} \frac{2\pi - T(\nu_j)}{T(\nu_j)} \theta(\nu_j, \D)\\
  & = \sum_{j=1}^q \frac{2 \pi - T(\nu_j)}{T(\nu_j)}\cdot \pi 
  = \sum_{j=1}^q \frac{2\pi -2\pi m_j}{2\pi m_j}\cdot\pi 
  = \pi \left( \sum_{j=1}^q \frac{1}{m_j} - q \right), 
 \end{aligned}
\end{equation}
where $v_i(\D)$'s are the vertices of $\D$. Here we used the convention that $\theta(v,\D) = 0$ if $v$ is not a vertex of $\D$,
and that $\infty / \infty = 1$ and $a / \infty = 0$ for $a < \infty$.

On the other hand, the Gaussian curvature of the pull-back spherical metric on $S$
exists and is 1 at non-singular points.
Hence if $D$ is a region in $S$ which contains no singular points, the 
integral curvature $\omega(D)$ of $D$ is same as the area of $D$. Since $H$
contains no singular point in its interior and is isometric to either
$G_1$ or $G_2$, we have
\begin{equation}\label{E:SumArea}
 \sum_{\D \subseteq H} \omega(\D^\circ) = \omega(H^\circ) = 2 \pi.
\end{equation}
Also note that the left turn of each side of triangles in $\T$ is zero because it is a 
locally geodesic curve containing no singular points on its interior arc. 
Therefore by the discussion in Section~\ref{theorem}, we have
\begin{equation}\label{mmmm}
 K(\D) = \sum_{i=1}^3 \frac{2\pi - T(v_i)}{T(v_i)} \theta(v_i,\D) + \omega(\D^\circ),
\end{equation} 
and by \eqref{E:SumCur}, \eqref{E:SumArea} and \eqref{mmmm}
we conclude that
\begin{align*}
 \sum_{\D \subseteq H} & K(\D)  = \sum_{\D \subseteq H} \left( \sum_{i=1}^3 
   \frac{2\pi - T(v_i(\D))}{T(v_i(\D))} \theta(v_i(\D), \D) + \omega(\D^\circ)
   \right) \\
 &  = \sum_{\D \subseteq H} \sum_{i=1}^3 \frac{2\pi - T(v_i(\D))}{T(v_i(\D))}
   \theta(v_i(\D), \D) + \sum_{\D \subseteq H} \omega(\D^\circ) \\
 & = \pi \left( \sum_{j=1}^q \frac{1}{m_j} - q \right) + 2 \pi
   = \pi \left( \sum_{j=1}^q \frac{1}{m_j} - q +2 \right) 
   = \pi E_p.
\end{align*}
Now it is easy to see that (M2)$'$ implies (M2). The theorem follows.
\end{proof}

This proof, as we discussed in Remark~\ref{rkk}, in fact shows that the surface $S \in F_q$ satisfies
a linear isoperimetric inequality with respect to the pull-back spherical metric.

\begin{remark}
Theorem~\ref{T2} can be regarded as a generalization of the hyperbolic case of 
Theorem~\ref{T:Ne}. First of all, Theorem~\ref{T:Ne} is nothing but the case
$M=1$ in Theorem~\ref{T2}. Secondly, $S$ does not have
to be regularly ramified in Theorem~\ref{T2} while it has to be in 
Theorem~\ref{T:Ne}. Also note that some points over $a_j$'s are 
even allowed to be unramified in Theorem~\ref{T2}. 
\end{remark}

We finish this paper with some discussion about Nevanlinna's conjecture. Suppose
a line complex $\Gamma$ of degree $q$ and the corresponding Riemann surface $S$ of class $F_q$ are given.
Let $V\Gamma$ be the vertex set of $\Gamma$ and
define the distance  $d_\Gamma(p,p')$ between two vertices $p, p' \in V\Gamma$ as the infimum of 
the (combinatorial) lengths of curves in $\Gamma$ connecting $p$ and $p'$. Now for
a vertex $p \in V\Gamma$, let
\[
B(p,j) := \{ p' \in V\Gamma : d_\Gamma(p, p') \leq j \}.
\]

R.~Nevanlinna defined \emph{the mean excess} of $\Gamma$ (or $S$) by the formula
\begin{equation}\label{meanexcess}
  E := \lim_{j \to \infty} \frac{1}{n_j} \sum_{p' \in B(p,j)} E_{p'},
\end{equation}
where $n_j$ is the number of vertices in $B(p,j)$. 
If this limit does not exist, one can only consider the upper ($\overline{E}$) or lower ($\underline{E}$)
excess of $S$, which are defined by \eqref{meanexcess} with limit superior or limit inferior
in place of limit, respectively.
Note that this definition of excess coincides with the quantity $E$ in
Theorem~\ref{T:Ne} when $S$ is regularly ramified.
In general $E$ depends on the base point $p$.

If $S$ is an $n$-sheeted closed surface, i.e., the corresponding line complex is a finite graph with
$2n$ vertices, 
one can show by Euler's formula that $E$ always exists and is 
equal to $2/n$. For an infinite-sheeted open surface $S$, Nevanlinna stated the following 
conjecture \cite{Ne}:
\emph{Is the surface parabolic if $E = 0$ and hyperbolic if $E<0$?}
Unfortunately, however, the answer is negative. Teichm\"uller \cite{Tei} 
constructed a hyperbolic surface with $E=0$ and recently 
I.~Benjamini, S.~Merenkov and O.~Schramm \cite{BMO} proved the existence of a parabolic surface with 
$\overline{E} < 0$. 

Therefore a stronger condition is necessary. Indeed, we prove the following
Theorem~\ref{T:final} with an additional ``uniform'' type assumption.

\begin{theorem}\label{T:final}
Suppose $\Gamma$ is a line complex and let $S$ be the associated Riemann
surface of class $F_q$. If there exist constants $\epsilon > 0$ and $M > 0$ such that
\[
 \sum_{p \in \Gamma_0} E_p \leq - \epsilon
\]
for every connected subgraph $\Gamma_0$ with $\# (\Gamma_0) \geq M$, then
$S$ is hyperbolic.
\end{theorem}

We need the following two lemmas for the proof of this theorem.

\begin{lemma}\label{L:Par}
Suppose $\Lambda$ is a connected subgraph of a line complex of degree $q$ 
such that $4q \leq \#(\Lambda) =: K < \infty$. Then there are two disjoint 
connected subgraphs $\Lambda_1$ and $\Lambda_2$ of $\Lambda$ such 
that $\Lambda_1 \cup \Lambda_2 = \Lambda$ and $\#(\Lambda_i) \geq K/2q$
for $i=1,2$.
\end{lemma}

\begin{proof}
We describe an algorithm to find such $\Lambda_1$ and $\Lambda_2$.

\textbf{Step 1}. Take a subgraph $\Lambda^1 \subseteq \Lambda$ such that
$\#(\Lambda^1) = 1$. If $\Lambda \backslash \Lambda^1$ is disconnected, we
set $\Lambda^0 = \emptyset$ and go to step 3 with $k=1$. Otherwise, we go to step 2 with $k=2$.

\textbf{Step 2}. Suppose that we have found a connected subgraph 
$\Lambda^{k-1}$ such that $\Lambda \backslash \Lambda^{k-1}$ is connected
and $\#(\Lambda^{k-1}) = k-1 < K/2q$. We then choose a connected subgraph
$\Lambda^k$ such that $\Lambda^{k-1} \subseteq \Lambda^k \subseteq \Lambda$
and $\#(\Lambda^k) = k$. If $\Lambda \backslash \Lambda^k$ is disconnected,
we go to step 3. If $\Lambda \backslash \Lambda^k$
is connected and $k = K/2q$, then the lemma follows with $\Lambda_1 = 
\Lambda^k$ and $\Lambda_2 = \Lambda \backslash \Lambda^k$. If
$\Lambda \backslash \Lambda^k$ is connected and $k < K/2q$, then we repeat
step 2.

\textbf{Step 3}. Let $p$ be the vertex of $\Lambda$ such that 
$\{p\} = \Lambda^k \backslash \Lambda^{k-1}$. Then since $\Lambda \backslash \Lambda^{k-1}$
is connected and there are at most $q$ edges at $p$, there are at
most $q$ components in $\Lambda \backslash \Lambda^k$. Moreover because $k < K/2q + 1 \leq K/2$, 
\[
\#(\Lambda \backslash \Lambda^k) = K - k \geq \frac{K}{2},
\]
hence there exists a component $\Lambda_1$ such that $\#(\Lambda_1) \geq
K/2q$. Let $\Lambda_2 := \Lambda \backslash \Lambda_1$. The subgraph $\Lambda_2$
is connected because each component in $\Lambda \backslash \Lambda^k$ is
connected by an edge to $p \in \Lambda^k$.
Now there are two cases. If $\#(\Lambda_2) \geq K/2q$, then $\Lambda_1$ and
$\Lambda_2$ satisfy all the properties in the lemma. If $\#(\Lambda_2) 
< K/2q$, let $\Lambda^{k'-1} := \Lambda_2$ where $k'-1 = \#(\Lambda_2)$,
and note that $k +1 \leq k'$ because $\Lambda^k \subseteq \Lambda_2$ by our construction. 
Now we repeat step 2 with $k' = k$. This is the end of all steps.

Note that this process must terminate since $K < \infty$. The lemma follows.
\end{proof}

\begin{lemma}\label{L:Par2}
Suppose $\Lambda$ is a connected subgraph of a line complex of degree $q$ such
that $\#(\Lambda) = \infty$. Then for any integer
$M \geq 2$, there are disjoint connected
subgraphs $\Lambda_1, \Lambda_2, \ldots, \Lambda_s$ for some $s$ such
that $\bigcup_{i=1}^s \Lambda_i = \Lambda$, $\#(\Lambda_1) < \infty$, and for 
each $i$, either $\#(\Lambda_i) = \infty$ or $M \leq \#(\Lambda_i) \leq 2qM^2$.
\end{lemma}

\begin{proof}
We first take a connected subgraph $\Lambda_1 \subseteq \Lambda$ such that
$\#(\Lambda_1) = M$. Since each component of $\Lambda \backslash \Lambda_1$
is connected by an edge to at least one vertex of $\Lambda_1$ and each vertex
has at most $q$ edges, $\Lambda \backslash \Lambda_1$ has 
at most $qM$ components. Let $\Lambda^1, 
\Lambda^2, \ldots, \Lambda^k$, $k \leq qM$, be the components of $\Lambda
\backslash \Lambda_1$. 

If $\#(\Lambda^j) = \infty$ or $M \leq \#(\Lambda^j) \leq 2qM$ for some 
$j$'s, we take each of these $\Lambda^j$'s as one of $\Lambda_i$, 
$i \geq 2$, in the lemma. If $2qM < \#(\Lambda^j) < \infty$ for some $j$'s, 
we apply Lemma~\ref{L:Par} to each of these and their subgraphs repeatedly 
to find disjoint connected 
subgraphs $\Lambda^{j1}, \Lambda^{j2}, \ldots, \Lambda^{jm_j}$ such that 
$M \leq \#(\Lambda^{jl}) \leq 2qM$, $l=1, \ldots, m_j$, and 
$\bigcup_{l=1}^{m_j} \Lambda^{jl} = \Lambda^j$. We then take each of 
these $\Lambda^{jl}$'s as one of $\Lambda_i$, $i \geq 2$. Finally, we replace
$\Lambda_1$ by $\Lambda_1 \cup \bigcup \Lambda^j$ where the union is over 
all $\Lambda^j$'s with $\#(\Lambda^j) < M$. Note that $\Lambda_1$ remains 
to be connected and $\#(\Lambda_1) \leq 2qM^2$ since $k$, the number of 
components in $\Lambda \backslash \Lambda_1$, is at most $qM$. This completes
the proof.
\end{proof}

\begin{proof}[Proof of Theorem~\ref{T:final}]
It is enough to find a partition $\PP(\Gamma)$ of $\Gamma$ satisfying 
\[
 M \leq  \#(\Gamma_\alpha) \leq 2qM^2
\]
for all $\Gamma_\alpha \in \PP(\Gamma)$. Conditions (M1)$'$ and (M2)$'$ 
in Theorem~\ref{T2} are trivially satisfied in this case.

Let $\mathcal{F} \subseteq 2^\Gamma$ be the set of all connected subgraphs
of $\Gamma$. Then we define 
\begin{align*}
 \mathfrak{P} := & \{ P \subseteq \mathcal{F} : \Gamma \backslash 
 \bigcup_{\Gamma_\alpha \in P} \Gamma_\alpha \mbox{ has no finite 
 component, and for all $\Gamma_\alpha$} \\
 &\mbox{ and $\Gamma_\beta$ in $P$, } 
 M \leq \#(\Gamma_\alpha) \leq 2qM^2 \mbox{ and } \Gamma_\alpha \cap 
 \Gamma_\beta = \emptyset \mbox{ if $\alpha \ne \beta$}\}. 
\end{align*}

The collection $\mathfrak{P}$ is nonempty since $\emptyset \in \mathfrak{P}$, and $\mathfrak{P}$
is partially ordered under usual set inclusion. Moreover, every chain
$P_1 \subseteq P_2 \subseteq \cdots$ in $\mathfrak{P}$ has an upper bound
$P'=\bigcup P_i$. In fact, if $\Lambda$ is a finite component of $\Gamma 
\backslash \bigcup_{\Gamma_\alpha \in P'} \Gamma_\alpha$, then it is 
connected by single edges to a finite number of vertices in 
$\bigcup_{\Gamma_\alpha \in P'} \Gamma_\alpha$. But all of these vertices 
must belong to $\bigcup_{\Gamma_\alpha \in P_n} \Gamma_\alpha$ for 
sufficiently large $n$, so $\Lambda$ is a finite component of $\Gamma 
\backslash \bigcup_{\Gamma_\alpha \in P_n} \Gamma_\alpha$. This is impossible
because $P_n \in \mathfrak{P}$, and one can easily see that $P' \in \mathfrak{P}$.

Now Zorn's lemma implies the existence of a maximal element $\PP(\Gamma)$
in $\mathfrak{P}$. To show that $\PP(\Gamma)$ is a desired partition,
it suffices to show $\bigcup_{\Gamma_\alpha \in \PP(\Gamma)} 
\Gamma_\alpha = \Gamma$. But if $\bigcup_{\Gamma_\alpha 
\in \PP(\Gamma)} \Gamma_\alpha \ne \Gamma$, then there is a component 
$\Lambda$ in $\Gamma \backslash \bigcup_{\Gamma_\alpha \in \PP(\Gamma)} 
\Gamma_\alpha$, which must be an infinite component by the definition of 
$\mathfrak{P}$. We then apply Lemma~\ref{L:Par2} to find a
partition $\{ \Lambda_1, \ldots, \Lambda_s \}$ of $\Lambda$ such that
$\#(\Lambda_1) < \infty$, and for each $i$, either $\#(\Lambda_i) = \infty$
or $M \leq \#(\Lambda_i) \leq 2qM^2$. Then
\[
 \PP(\Gamma) \cup \{\Lambda_i : M \leq \#(\Lambda_i) \leq 2qM^2 \}
\]
is an element in $\mathfrak{P}$ which is strictly larger than $\PP(\Gamma)$
because $\#(\Lambda_1) < \infty$. Hence $\PP(\Gamma)$ is not a maximal
element, which is a contradiction. The theorem follows.
\end{proof}

In the proof of Theorem~\ref{T:final} we only used the fact that $\Gamma$
is a connected graph of degree $q$. But each triangulation of a surface
$S$ can be represented by a graph of degree 3, 
so the following theorem may be proved by
using the same method.

\begin{theorem}
Suppose that $S$ is an open simply-connected Aleksandrov surface of curvature
at most $k$. If $S$ has a tiling
$\T$ satisfying the conditions \emph{(R1)} and \emph{(R2)} in 
Theorem~\ref{T} and there exist constants $\epsilon > 0$ and $M > 0$ such
that
\[
 \sum_{\D \subseteq C_\alpha} K(\D) \leq - \epsilon
\]
for every connected cluster $C_\alpha$ with $\# (C_\alpha) \geq M$, then
$S$ is hyperbolic.
\end{theorem}


\end{document}